\documentclass[12pt]{article}

\usepackage{amsmath, epsfig, cite,  setspace}
\usepackage{amssymb}
\usepackage{amsfonts}
\usepackage{latexsym}
\usepackage{amsthm}

\makeatletter
\renewcommand{\@seccntformat}[1]{{\csname the#1\endcsname}.\hspace{.5em}}
\makeatother

\newtheorem{thm}{Theorem}[section]

\newtheorem{cor}[thm]{Corollary}
\newtheorem{conj}[thm]{Conjecture}
\newtheorem{lem}[thm]{Lemma}

\renewcommand{\thefootnote}{*}

\numberwithin{equation}{section}

\begin{document}

\begin{center}
{\large\bf Further $q$-supercongruences from a transformation of Rahman }
\end{center}

\vskip 2mm \centerline{Victor J. W. Guo}
\begin{center}
{\footnotesize School of Mathematics and Statistics, Huaiyin Normal
University, Huai'an 223300, Jiangsu, People's Republic of China\\
{\tt jwguo@hytc.edu.cn }  }
\end{center}

\vskip 0.7cm \noindent{\bf Abstract.} Employing a quadratic  transformation formula of Rahman and the method of `creative microscoping'
(introduced by the author and Zudilin in 2019), we provide some new $q$-supercongruences for truncated basic hypergeometric
series. In particular, we confirm two recent conjectures of Liu and Wang.
We also propose some related conjectures on supercongruences and $q$-supercongruences.

\vskip 3mm \noindent {\it Keywords}: cyclotomic polynomials;  $q$-supercongruences; Rahman's transformation; creative microscoping

\vskip 0.2cm \noindent{\it AMS Subject Classifications}: 33D15; Secondary 11A07, 11B65

\renewcommand{\thefootnote}{**}

\section{Introduction}

In his first letter to Hardy on the 16th January 1913,  Ramanujan mentioned the following formula:
\begin{align}
\sum_{k=0}^\infty(8k+1)\frac{(\frac{1}{4})_k^4}{k!^4}
=\frac{2\sqrt{2}}{\sqrt{\pi}\,\Gamma(\frac 34)^2} \label{eq:ram1}
\end{align}
(see \cite[p.~25, (2)]{BR}), where $(a)_n=a(a+1)\cdots(a+n-1)$ is the Pochhammer symbol and $\Gamma(x)$ denotes the Gamma function.
Recently, Chen and Chu\cite{Chu1} gave a $q$-analogue of \eqref{eq:ram1} as follows:
\begin{equation}
\sum_{k=0}^{\infty}[8k+1] \frac{(q ; q^{4})_{k}^4}{(q^4 ; q^{4})_k^4}q^{2k}
=\frac{(q^5,q^3,q^3,q^3;q^4)_\infty}{(q^2,q^4,q^4,q^4;q^4)_\infty}.\label{eq:Chu1}
\end{equation}
They also obtained the following similar formula:
\begin{align}
\sum_{k=0}^{\infty} [6k+1] \frac{(q ; q^{4})_{k}(q ; q^{2})_{k}^{3}}{(q^{2} ; q^{2})_{k}(q^{4} ; q^{4})_{k}^{3}} q^{k^{2}+k}
=\frac{(q^5,q^3,q^3,q^3;q^4)_\infty}{(q^2,q^4,q^4,q^4;q^4)_\infty}.\label{eq;Chu2}
\end{align}
Here and in what follows, $(a;q)_n=(1-a)(1-aq)\cdots (1-aq^{n-1})$
is the {\em $q$-shifted factorial}, and $[n]:=[n]_q=(1-q^n)/(1-q)$ denotes the {\em $q$-integer}.
For simplicity, we also write $(a_1,a_2,\ldots,a_m;q)_n=(a_1;q)_n (a_2;q)_n\cdots (a_m;q)_n$ for
$n\geqslant 0 $ or $n=\infty$.

For any odd prime $p$, let $\Gamma_p(x)$ be the $p$-adic Gamma function \cite{Morita}. In 2015, Swisher \cite{Swisher} proved the following $p$-adic analogue of \eqref{eq:ram1}:
\begin{align}
\sum_{k=0}^{(p-1)/4}(8k+1)\frac{(\frac{1}{4})_k^4}{k!^4}
\equiv p\frac{\Gamma_p(\frac 12)\Gamma_p(\frac 14)}{\Gamma_p(\frac 34)}
\pmod{p^3},\quad\text{if $p\equiv 1\pmod{4}$},
\label{eq:hamme}
\end{align}
which was originally conjectured by Van Hamme \cite[(G.2)]{Hamme}.
Liu and Wang \cite{LiuW2} showed that \eqref{eq:hamme} can also be deduced from the following $q$-supercongruence:
\begin{align}
\sum_{k=0}^{(n-1)/4}[8k+1] \frac{(q ; q^{4})_{k}^4}{(q^4 ; q^{4})_k^4}q^{2k}
\equiv \frac{(q^{2} ; q^{4})_{(n-1)/{4}}}{(q^{4} ; q^{4})_{(n-1)/4}}[n] q^{(1-n) / 4}  \pmod {[n]\Phi_{n}(q)^{2}}\label{eq;qG2}
\end{align}
(for a more general form, see \cite[Theorem 4.3]{GS0}). On the other hand, Guo and Schlosser  \cite[Theorem 2 with $d=4$]{GS2} proved that,
\begin{align*}
\sum_{k=0}^{n-1}[8k+1]\frac{(q;q^4)_k^4}{(q^4;q^4)_k^4}q^{2k}
\equiv 0\pmod{\Phi_n(q)^2} \quad\text{for $n\equiv 3 \pmod 4$.}
\end{align*}
Here, $\Phi_n(q)$ is the $n$-th {\em cyclotomic polynomial} in $q$, which may be given by
\begin{align*}
\Phi_n(q)=\prod_{\substack{1\leqslant k\leqslant n\\ \gcd(n,k)=1}}(q-\zeta^k),
\end{align*}
where $\zeta$ is an $n$-th primitive root of unity. Moreover, two rational functions $A(q)$ and $B(q)$ in $q$ are called congruent modulo a polynomial
$P(q)$, denoted by $A(q) \equiv B(q) \pmod{P(q)}$, if $P(q)$ divides the numerator of the reduced form
of $A(q)-B(q)$ in the polynomial ring $\mathbb{Z}[q]$.

Recently, Liu and Wang \cite{LiuW3} proved that,
modulo $[n]\Phi_{n}(q)^2$,
\begin{align}
\sum_{k=0}^{M}[6k+1] \frac{\left(q ; q^{4}\right)_{k}\left(q ; q^{2}\right)_{k}^{3}}{\left(q^{2} ; q^{2}\right)_{k}\left(q^{4} ; q^{4}\right)_{k}^{3}} q^{k^{2}+k}
\equiv \begin{cases}
\dfrac{\left(q^{2} ; q^{4}\right)_{(n-1)/4}}{\left(q^{4} ; q^{4}\right)_{(n-1)/4}}[n] q^{(1-n) / 4} ,&\text{if $n\equiv 1 \pmod{4}$,} \\
0 ,&\text{if $n\equiv 3 \pmod{4}$},
\end{cases}\label{eq:6k+1-LW}
\end{align}
where $M=(n-1)/2$ or $n-1$.
They also gave the following generalization of the second case of  \eqref{eq:6k+1-LW} modulo $\Phi_n(q)^2$,
for any positive integer $d\geqslant 2$ and positive odd integer $n$ with $ n\equiv d+1 \pmod {2d}$,
\begin{align}
\sum_{k=0}^{(n-1)/d}[3dk+1] \frac{(q,e,q^{1+d}/e; q^{2d})_{k}(q,q,q^{d-1};q^{d})_{k}}
{(q^{d},e,q^{1+d}/e;q^{d})_{k}(q^{2d},q^{2d},q^{d+2};q^{2d})_{k}}q^{dk}\equiv 0 \pmod{\Phi_{n}(q)^2}. \label{eq:LW}
\end{align}
For some other recent work on $q$-supercongruences, see
\cite{Mohamed,Guo-aam,GL,GS3,GuoZu,GuoZu2,Liu,LP,LiuW,NW,SW,Wang,Wei,Zu19}.

The first aim of this paper is to establish the a stronger version of \eqref{eq:LW} for even $d$ as follows.
\begin{thm}\label{thm:main}
Let $d\geqslant 2$ be an even integer and $e$ an indeterminate. Let $n\equiv d+1\pmod {2d}$ be a positive integer. Then
\begin{align}
\sum_{k=0}^{(dn+n-1)/(2d)}
[3dk+1] \frac{(q,e,q^{1+d}/e; q^{2d})_{k}(q,q,q^{d-1};q^{d})_{k}}
{(q^{d},e,q^{1+d}/e;q^{d})_{k}(q^{2d},q^{2d},q^{d+2};q^{2d})_{k}}q^{dk}\equiv 0 \pmod{\Phi_{n}(q)^3}.  \label{eq:main}
\end{align}
\end{thm}

Note that the $q$-supercongruence \eqref{eq:main} modulo $\Phi_n(q)^2$ also follows from \eqref{eq:LW}, since
the $k$-th summand on the left-hand side of \eqref{eq:LW} is congruent to $0$ modulo $\Phi_n(q)^2$ for all $(n-1)/d<k\leqslant (dn+n-1)/(2d)$.
However, the $q$-supercongruence \eqref{eq:LW} does not hold modulo $\Phi_n(q)^3$ in general.

Letting $e\to 0$ and $e=-q$ in \eqref{eq:main}, respectively, we obtain
\begin{align*}
\sum_{k=0}^{(dn+n-1)/(2d)}
[3dk+1] \frac{(q; q^{2d})_{k}(q,q,q^{d-1};q^{d})_{k}}
{(q^{d};q^{d})_{k}(q^{2d},q^{2d},q^{d+2};q^{2d})_{k}}q^{d(k^2+k)/2}\equiv 0 \pmod{\Phi_{n}(q)^3}, \\
\sum_{k=0}^{(dn+n-1)/(2d)}
[3dk+1] \frac{(q,-q,-q^d; q^{2d})_{k}(q,q,q^{d-1};q^{d})_{k}}
{(q^{d},-q,-q^d;q^{d})_{k}(q^{2d},q^{2d},q^{d+2};q^{2d})_{k}}q^{dk}\equiv 0 \pmod{\Phi_{n}(q)^3}.
\end{align*}

Letting $n=p$ be an odd prime and $q\to  1$ in each of the above $q$-supercongruences, we get the following result: for even $d\geqslant 2$ and $p\equiv d+1\pmod{2d}$,
\begin{align}
\sum_{k=0}^{(dp+p-1)/(2d)} (3dk+1)\frac{(\frac{1}{d})_{k}^{2}(\frac{d-1}{d})_{k}(\frac{1}{2d})_{k}}
{k!^{3}4^{k} (\frac{d+2}{2d})_{k}}
&\equiv 0 \pmod {p^3}. \label{eq:dp+p-1}
\end{align}
For $d=4$, we have $\frac{d-1}{d}=\frac{d+2}{2d}=\frac{3}{4}$. It is easy to see $(\frac{1}{2})_k\equiv (\frac{1}{4})_k\equiv 0\pmod{p}$ for
$(3p-1)/4<k\leqslant p-1$. Thus, from the $d=4$ case of \eqref{eq:dp+p-1} we deduce the following supercongruence, which was conjectured by Liu and Wang \cite[Conjecture 4]{LiuW3}.
\begin{cor}
Let $p \equiv 5 \pmod 8 $ be a prime. Then
\begin{align*}
 \sum_{k=0}^{p-1}(12 k+1) \frac{(\frac{1}{4})_{k}^{2}(\frac{1}{8})_{k}}{ k!^3 4^{k}} \equiv 0 \pmod {p^3}.
\end{align*}
\end{cor}

The second aim of this paper is to establish the following $q$-supercongruence, which is a generalization of \cite[Theorm 8]{LiuW3} for the second case with $d$ even.
\begin{thm}\label{thm:main2}
Let $d\geqslant 2$ be an even integer and $e$ an indeterminate. Let $n\equiv d+1\pmod{2d}$ be a positive integer. Then, modulo $\Phi_{n}(q)^3$,
\begin{align}
\sum_{k=0}^{M}[3dk-1]\frac{(q^{-1},e,q^{d-1}/e; q^{2d})_{k}(q,q,q^{d-3};q^{d})_{k} q^{dk} }
{(q^{d},e,q^{d-1}/e;q^{d})_{k}(q^{2d-2}, q^{2d-2},q^{d+2};q^{2 d})_{k} }
\equiv 0,   \label{eq:main2}
\end{align}
where $M=(n-1)/2$ if $d=2$, and $M=(dn-n+1)/(2d)$ otherwise.
\end{thm}

The third aim of this paper is to prove the following $q$-supercongruence, which was originally conjectured by Liu and Wang \cite[Conjecture 5]{LiuW3}.
\begin{thm}\label{thm:main3}
Let $n$ be a positive odd integer. Then, modulo $\Phi_{n}(q)^3$,
\begin{align}
&\sum_{k=0}^{(n-1)/2}[6k-1]\frac{(q^{-1};q^{4})_k(q^{-1};q^2)_k(q;q^2)_k^2}{(q^2;q^2)_k(q^{4};q^{4})_k(q^{2};q^{4})_k^2}q^{k^2+k+1} \notag \\
&\quad\equiv \begin{cases}
-\dfrac{(q^{2} ; q^{4})_{(n-1)/4}}{(q^{4} ; q^{4})_{(n-1)/4}} q^{(n-1) / 4} ,& \text{if $n\equiv 1 \pmod{4}$,} \\[10pt]
0,& \text{if $n\equiv 3 \pmod{4}$}.
\end{cases}  \label{conj:cases}
\end{align}
\end{thm}
Liu and Wang originally conjectured that \eqref{conj:cases} holds modulo $[n]\Phi_n(q)^2$, which is not true (the first counterexample is $n=15$).

We shall prove Theorems \ref{thm:main}, \ref{thm:main2}, and \ref{thm:main3} by using the method of `creative microscoping', which was introduced by the author and Zudilin \cite{GuoZu}.
At the end of this paper, we put forward several open problems on supercongruences and $q$-supercongruences.

\section{Proof of Theorem \ref{thm:main} }  \label{sec:2}

Recall that the {\it basic hypergeometric series $_{r+1}\phi_r$} (see Gasper and Rahman's monograph \cite{GR}) is defined as
$$
_{r+1}\phi_{r}\!\left[\begin{array}{c}
a_1,a_2,\ldots,a_{r+1}\\
b_1,b_2,\ldots,b_{r}
\end{array};q,\, z
\right]
=\sum_{k=0}^{\infty}\frac{(a_1,a_2,\ldots, a_{r+1};q)_k z^k}
{(q,b_1,\ldots,b_{r};q)_k}.
$$
We need a quadratic transformation of Rahman \cite[(3.8.13)]{GR}, which can be stated as follows:
\begin{align}
&\sum_{k=0}^{\infty} \frac{(1-a q^{3 k})(a, d, a q / d ; q^{2})_{k}(b, c, a q / b c ; q)_{k}}
{(1-a)(a q / d, d, q ; q)_{k}(a q^{2} / b, a q^{2} / c, b c q ; q^{2})_{k}} q^{k}\nonumber \\
&=\frac{(a q^{2}, b q, c q, a q^{2} / b c ; q^{2})_{\infty}}{(q, a q^{2} / b, a q^{2} / c, b c q ; q^{2})_{\infty}}
\,{ }_{3} \phi_{2}\!\left[\begin{array}{c}
b, c, a q / b c \\
d q, a q^{2} / d
\end{array} ; q^{2}, q^{2}\right],
\label{eq;Rm}
\end{align}
provided that $d$ or $aq/d$ is not of the form $q^{-2n}$ ($n$ is a non-negative integer).

We first give a generalization of Theorem \ref{thm:main} with an extra parameter $a$. Note that this $q$-congruence modulo $(1-aq^n)(a-q^n)$ was
already indicated by Liu and Wang \cite{LiuW3}. In order to make the paper self-contained, we give a complete proof here.
\begin{thm}\label{lem:second}
Let $d\geqslant 2$ be an even integer and $a,\,e$ indeterminates. Let $n\equiv d+1\pmod {2d}$ be a positive integer. Then, modulo $\Phi_n(q)(1-aq^n)(a-q^n)$,
\begin{equation}
\sum_{k=0}^{(dn+n-1)/(2d)}[3d k+1] \frac{(q, e,  q^{1+d} / e ; q^{2 d})_{k}(a q, q / a, q^{d-1} ; q^{d})_{k}}
{(q^{d}, e,  q^{1+d} / e ; q^{d})_{k}(a q^{2 d}, q^{2 d} / a, q^{d+2} ; q^{2 d})_{k}} q^{dk}
\equiv 0\label{eq:lemma2}.
\end{equation}
\end{thm}

\begin{proof}
Put $q \mapsto q^{d}$, $a=q^{1-dn-n}$,  $b=aq$, $c=q/a$ and $d=e$ in Rahman's transformation \eqref{eq;Rm}. Then, for $n\equiv d+1\pmod {2d}$, we have
\begin{align}
&\sum_{k=0}^{(dn+n-1)/(2d)}\frac{(1-q^{3dk+1-dn-n})(q^{1-dn-n},e,q^{1+d-dn-n}/e;q^{2d})_k(aq,q/a,q^{d-1-dn-n};q^d)_k}
{(1-q^{1-dn-n})(q^{d},e,q^{1+d-dn-n}/e ; q^{d})_k(aq^{2d-dn-n},q^{2d-dn-n}/a,q^{d+2};q^{2d})_k}q^{dk} \notag \\
&=\frac{(q^{2d+1-dn-n},q^{2d-1-dn-n},aq^{d+1},q^{d+1}/a;q^{2d})_\infty}{(q^d,q^{d+2},aq^{2d-dn-n},q^{2d-dn-n}/a;q^{2d})_\infty}
\sum_{k=0}^{\infty}\frac{(aq,q/a,q^{d-1-dn-n};q^{2d})_k}{(q^{2d},eq^{d},q^{2d+1-dn-n}/e;q^{2d})_k}q^{2dk} \notag\\
&=0,  \label{eq:special}
\end{align}
where we have used  $(q^{2d+1-dn-n};q^{2d})_\infty=0$ for $n\equiv d+1\pmod{2d}$. It is easy to see that, for $0\leqslant k\leqslant (dn+n-1)/(2d)$,
the polynomial
$$
(1-q^{1-dn-n})(q^{d},e,q^{1+d-dn-n}/e ; q^{d})_k(aq^{2d-dn-n},q^{2d-dn-n}/a,q^{d+2};q^{2d})_k
$$
is relatively prime to $\Phi_n(q)$.
Since $q^{n} \equiv 1\pmod {\Phi_{n}(q)}$, from \eqref{eq:special} we deduce that the $q$-congruence \eqref{eq:lemma2} is true modulo $\Phi_n(q)$.

On the other hand, letting $q \mapsto q^{d}, a=q, b=q^{1-n}$, $c=q^{1+n}$ and $d=e$ in \eqref{eq;Rm}, we get
\begin{align*}
&\sum_{k=0}^{(n-1)/d}[3d k+1] \frac{(q, e,  q^{1+d} / e ; q^{2 d})_{k}(q^{1-n}, q^{1+n}, q^{d-1} ; q^{d})_{k}}
{(q^{d}, e,  q^{1+d} / e ; q^{d})_{k}(q^{2d-n}, q^{2d+n}, q^{d+2} ; q^{2 d})_{k}} q^{dk} \\
&=\frac{(q^{2d+1},q^{d+1-n},q^{d+1+n},q^{2d-1};q^{2d})_\infty}{(q^d,q^{2d+n},q^{2d-n},q^{d+2};q^{2d})_\infty}
\sum_{k=0}^{\infty}\frac{(q^{1-n},q^{1+n},q^{d-1};q^{2d})_k}{(q^{2d},eq^{d},q^{2d+1}/e;q^{2d})_k}q^{2dk}\\
&=0,
\end{align*}
since $(q^{d+1-n};q^{2d})_\infty=0$ for  $n\equiv d+1 \pmod {2d}$. Noticing that $(dn-1)/(2d)> (n-1)/d$,
we conclude that the left-hand side of \eqref{eq:lemma2}
is equal to $0$ for $a=q^{-n}$ and $a=q^n$. Namely, the $q$-congruence \eqref{eq:lemma2} is true modulo $1-aq^n$ and $a-q^n$.
Since $\Phi_n(q)$, $1-aq^n$, and $a-q^n$ are pairwise relatively prime polynomials in $q$, we complete the proof.
\end{proof}

\begin{proof}[Proof of Theorem {\rm\ref{thm:main}}]
Since $n\equiv d+1\pmod{2d}$, we have $\gcd(2d,n)=1$. Hence, $(q^{2d};q^{2d})_k$ is relatively prime to $\Phi_n(q)$ for any $0\leqslant k\leqslant n-1$.
It is clear that $(dn+n-1)/(2d)\leqslant n-1$. Moreover, the polynomial $1-q^n$ contains the factor $\Phi_n(q)$. The proof of \eqref{eq:main} then follows
from \eqref{eq:lemma2} by taking $a=1$.
\end{proof}

\section{Proof of Theorem \ref{thm:main2} }  \label{sec:3}

Like before, we first establish the following generalization of Theorem \ref{thm:main2} with an additional parameter $a$.
\begin{thm}\label{thm:main2-a}
Let $d\geqslant 2$ be an even integer and $a,e$ indeterminates. Let $n\equiv d+1\pmod{2d}$ be a positive integer. Then, modulo $\Phi_{n}(q)(1-aq^n)(a-q^n)$,
\begin{align}
\sum_{k=0}^{M}[3dk-1]\frac{(q^{-1},e,q^{d-1}/e; q^{2d})_{k}(aq,q/a,q^{d-3};q^{d})_{k} q^{dk} }
{(q^{d},e,q^{d-1}/e;q^{d})_{k}(aq^{2d-2}, q^{2d-2}/a,q^{d+2};q^{2 d})_{k} }
\equiv 0,   \label{eq:main2-a}
\end{align}
where $M=(n-1)/2$ if $d=2$, and $M=(dn-n+1)/(2d)$ otherwise.
\end{thm}

\begin{proof}
Set $q \mapsto q^{d}$, $a=q^{-1-dn+n}$,  $b=aq$, $c=q/a$ and $d=e$ in \eqref{eq;Rm}. Then, for $n\equiv d+1\pmod {2d}$, we have
\begin{align}
&\sum_{k=0}^{(dn-n+1)/(2d)}\frac{(1-q^{3dk-1-dn+n})(q^{-1-dn+n},e,q^{d-1-dn+n}/e;q^{2d})_k(aq,q/a,q^{d-3-dn+n};q^d)_k}
{(1-q^{-1-dn+n})(q^{d},e,q^{d-1-dn+n}/e; q^{d})_k(aq^{2d-2-dn+n},q^{2d-2-dn+n}/a,q^{d+2};q^{2d})_k}q^{dk} \notag \\
&=\frac{(q^{2d-1-dn+n},q^{2d-3-dn+n},aq^{d+1},q^{d+1}/a;q^{2d})_\infty}{(q^d,q^{d+2},aq^{2d-dn+n},q^{2d-dn+n}/a;q^{2d})_\infty}
\sum_{k=0}^{\infty}\frac{(aq,q/a,q^{d-3-dn+n};q^{2d})_k}{(q^{2d},eq^{d},q^{2d-1-dn+n}/e;q^{2d})_k}q^{2dk} \notag\\
&=0,  \label{eq:special-2}
\end{align}
where we have utilized  $(q^{2d-1-dn+n};q^{2d})_\infty=0$ for $n\equiv d+1\pmod{2d}$. Moreover, it is not difficult to
see that $(q^{2d-2};q^{2d})_k$ is relatively prime to $\Phi_n(q)$ for $0\leqslant k\leqslant M$ (in fact, this is true for $0\leqslant k\leqslant (d-1)(n-1)/d$).
Noticing $q^n\equiv 1\pmod{\Phi_n(q)}$ again, the modulus $\Phi_n(q)$ case of the $q$-congruence \eqref{eq:main2-a}  follows from \eqref{eq:special-2} immediately
(for $d=2$, we need to use the fact that $(q^{-1};q^4)_k\equiv 0\pmod{\Phi_n(q)}$ for $(n+1)/4<k\leqslant (n-1)/2$).

On the other hand, letting $q \mapsto q^{d}, a=q^{-1}, b=q^{1-n}$, $c=q^{1+n}$ and $d=e$ in \eqref{eq;Rm}, we obtain
\begin{align*}
&\sum_{k=0}^{(n-1)/d}[3dk-1]\frac{(q^{-1},e,q^{d-1}/e;q^{2d})_{k}(q^{1-n}, q^{1+n}, q^{d-3}; q^{d})_{k}}
{(q^{d},e,q^{d-1}/e;q^{d})_{k}(q^{2d-2-n}, q^{2d-2+n},q^{d+2};q^{2d})_{k}} q^{dk}\\[5pt]
&\quad=-\frac{(q^{2d-1},q^{2d-3},q^{d+1-n},q^{d+1+n};q^{2d})_\infty}{q(q^d,q^{d+2},q^{2d-2-n},q^{2d-2+n};q^{2d})_\infty}
\sum_{k=0}^{(n-1) /(2 d)} \frac{(q^{1-n}, q^{1+n}, q^{d-3} ; q^{2 d})_{k}}{(q^{2 d}, e q^{d}, q^{2 d-1}/e;q^{2d})_{k}} q^{2dk}\\[5pt]
&\quad=0,
\end{align*}
as was first given by Liu and Wang \cite{LiuW3}. In view of $(dn-n-1)/(2d)< (n-1)/2$ for $d=2$, and $(dn-n-1)/(2d)> (n-1)/d$ for $d\geqslant 4$, one sees that the left-hand side of \eqref{eq:main2-a}
is equal to $0$ for $a=q^{-n}$ and $a=q^n$. Thus, the $q$-congruence \eqref{eq:main2-a} holds modulo $1-aq^n$ and $a-q^n$.
Since the polynomials $\Phi_n(q)$, $1-aq^n$, and $a-q^n$ are relatively prime to one another, we accomplish the proof.
\end{proof}

\begin{proof}[Proof of Theorem {\rm\ref{thm:main2}}]
In the proof of Theorem \ref{thm:main2-a}, we have mentioned that $(q^{2d-2};q^{2d})_k$ is relatively prime to $\Phi_n(q)$ for $0\leqslant k\leqslant (dn-n+1)/(2d)$.
The proof of \eqref{eq:main2} then follows from the $a=1$ case of \eqref{eq:main2-a}.
\end{proof}

\section{Proof of Theorem \ref{thm:main3} }  \label{sec:4}

We require the following result, which was first given in \cite[Lemma 2.1]{Guo-aam}.
For the reader's convenience, we include a short proof here.
\begin{lem}\label{eq:lem:old}
Let $n$ be a positive odd integer and $a$ an indeterminate. Then
\begin{align}
(aq,q/a;q^2)_{(n-1)/2}
\equiv  (-1)^{(n-1)/2}\frac{(1-a^n)q^{(1-n^2)/4}}{(1-a)a^{(n-1)/2}} \pmod{\Phi_n(q)}.  \label{eq:lem-2}
\end{align}
\end{lem}
\begin{proof}It is easy to see that
\begin{align*}
(q/a;q^2)_{(n-1)/2}
&=(1-q/a)(1-q^3/a)\cdots(1-q^{n-2}/a) \\
&\equiv (1-q^{1-n}/a)(1-q^{3-n}/a)\cdots (1-q^{-2}/a)\\
&=(-1)^{(n-1)/2}(aq^2;q^2)_{(n-1)/2}\frac{q^{(1-n^2)/4}}{a^{(n-1)/2}} \pmod{\Phi_n(q)}.
\end{align*}
Hence, the left-hand side of \eqref{eq:lem-2} is congruent to
\begin{align*}
(-1)^{(n-1)/2}(aq;q)_{n-1}\frac{q^{(1-n^2)/4}}{a^{(n-1)/2}}.
\end{align*}
For any $n$-th primitive root of unity $\zeta$, we have
$$
(a\zeta;\zeta)_{n-1}=\frac{(a;\zeta)_n}{1-a}=\frac{1-a^n}{1-a},
$$
and so $(aq;q)_{n-1}$ is congruent to $(1-a^n)/(1-a)$ modulo $\Phi_n(q)$. This completes the proof.
\end{proof}

We have the following parametric generalization of Theorem \ref{thm:main3} for $n\equiv 1 \pmod{4}$.

\begin{thm}\label{lem:d=2}
Let $n\equiv 1\pmod{4}$ be a positive integer and $a$ an indeterminate. Then, modulo $\Phi_{n}(q)(1-aq^n)(a-q^n)$,
\begin{align}
\sum_{k=0}^{(n-1)/2}[6k-1]\frac{(q^{-1};q^{4})_k(q^{-1},aq,q/a;q^2)_k}
{(q^2;q^2)_k(q^{4},aq^2,q^2/a;q^{4})_k}q^{k^2+k+1}
\equiv \frac{(q^{2};q^{4})_{(n-1)/4}}{(q^{4};q^{4})_{(n-1)/4}} q^{(n-1)/4}.  \label{eq:d=2}
\end{align}
\end{thm}

\begin{proof}
Letting $d\to 0$ in \eqref{eq;Rm}, we have
\begin{align}
\sum_{k=0}^{\infty} \frac{(1-a q^{3 k})(a;q^{2})_{k}(b,c,aq/bc; q)_{k}}
{(1-a)(q;q)_{k}(a q^{2}/b,aq^{2}/c, bcq;q^{2})_{k}}q^{(k^2+k)/2}
=\frac{(aq^{2}, bq, cq, aq^{2}/bc; q^{2})_{\infty}}{(q, aq^{2}/b, aq^{2}/c,bcq; q^{2})_{\infty}}. \label{eqd=0}
\end{align}
We then take $q \mapsto q^2$, $a=q^{-1-n}$,  $b=aq$, $c=q/a$ in the above formula to obtain
\begin{align}
&\sum_{k=0}^{(n+1)/2}\frac{(1-q^{6k-1-n})(q^{-1-n};q^4)_k(aq,q/a,q^{-1-n};q^2)_k}
{(1-q^{-1-n})(q^2; q^2)_k(aq^{2-n},q^{2-n}/a,q^4;q^4)_k}q^{k^2+k} \notag \\
&=\frac{(q^{3-n},aq^3,q^3/a,q^{1-n};q^4)_\infty}{(q^2,q^{2-n}/a,aq^{2-n},q^4;q^4)_\infty}\notag\\
&=0.
\end{align}
Since $q^n\equiv 1\pmod{\Phi_n(q)}$, we conclude from the above equality that
\begin{align*}
\sum_{k=0}^{(n+1)/2}[6k-1]\frac{(q^{-1};q^{4})_k(q^{-1},aq,q/a;q^2)_k}
{(q^2;q^2)_k(q^{4},aq^2,q^2/a;q^{4})_k}q^{k^2+k+1}
\equiv 0 \pmod{\Phi_n(q)}.
\end{align*}
Namely,
\begin{align}
&\sum_{k=0}^{(n-1)/2}[6k-1]\frac{(q^{-1};q^{4})_k(q^{-1},aq,q/a;q^2)_k}
{(q^2;q^2)_k(q^{4},aq^2,q^2/a;q^{4})_k}q^{k^2+k+1}  \notag\\
&\quad \equiv - [2]\frac{(q^{-1};q^{4})_{(n+1)/2}(q^{-1},aq,q/a;q^2)_{(n+1)/2}}
{(q^2;q^2)_{(n+1)/2}(q^{4},aq^2,q^2/a;q^{4})_{(n+1)/2}}q^{(n+1)(n+3)/4+1} \pmod{\Phi_n(q)}. \label{eq:negative}
\end{align}

By Lemma \ref{eq:lem:old}, we have
\begin{align*}
\frac{(aq,q/a;q^2)_{(n+1)/2}}{(aq^2,q^2/a;q^4)_{(n+1)/2}}
&=\frac{(aq,q/a;q^2)_{(n-1)/2}(1-aq^n)(1-q^n/a)}{(aq^2,q^2/a;q^4)_{(n-1)/2}(1-aq^{2n})(1-q^{2n}/a)}\\
&=q^{(n^2-1)/4} \pmod{\Phi_n(q)}.
\end{align*}
Moreover, modulo $\Phi_n(q)$,
\begin{align*}
\frac{(q^{-1};q^2)_{(n+1)/2}} {(q^2;q^2)_{(n+1)/2}}
&=\frac{(1-q^{-1})(1-q)\cdots (1-q^{n-2})}{(1-q^2)(1-q^4)\cdots(1-q^{n+1})}
\equiv (-1)^{(n+1)/2}q^{-(n+1)(n+3)/4}, \\
-[2]\frac{(q^{-1};q^4)_{(n+1)/2}}{(q^4;q^4)_{(n+1)/2}}q
&\equiv\frac{(q^{3};q^4)_{(n-1)/2}} {(q^4;q^4)_{(n-1)/2}}
=\frac{(q^{3};q^4)_{(n-1)/4}(q^{n+2};q^4)_{(n-1)/4}} {(q^4;q^4)_{(n-1)/4}(q^{n+3};q^4)_{(n-1)/4}}\\
&\equiv \frac{(q^{2};q^{4})_{(n-1)/4}}{(q^{4};q^{4})_{(n-1)/4}} .
\end{align*}
Employing the above three $q$-congruences, we see that the right-hand side of \eqref{eq:negative} reduces to
\begin{align*}
\frac{(q^{2};q^{4})_{(n-1)/4}}{(q^{4};q^{4})_{(n-1)/4}} q^{(n-1)/4} \pmod{\Phi_n(q)}.
\end{align*}
This proves that \eqref{eq:d=2} is true modulo $\Phi_n(q)$.

The modulus $(1-aq^n)(a-q^n)$ case of \eqref{eq:d=2} was already given by Liu and Wang \cite[(4.2) with $e\to 0$]{LiuW3}, and
this can be easily checked by putting $q \mapsto q^2$, $a=q^{-1}$,  $b=q^{1-n}$, $c=q^{1+n}$ in \eqref{eqd=0}.
Since the polynomials $\Phi_n(q)$ and $(1-aq^n)(a-q^n)$ are relatively prime, we finish the proof of the theorem.
\end{proof}

\begin{proof}[Proof of Theorem {\rm\ref{thm:main3}}]
Letting $a=1$ in \eqref{eq:d=2}, we arrive at \eqref{conj:cases} for $n\equiv 1\pmod{4}$.
On the other hand, letting $d=2$ and $e\to $ in \eqref{eq:main2}, we are led to \eqref{conj:cases} for $n\equiv 3\pmod{4}$.
\end{proof}

\section{Concluding remarks and open problems}  \label{sec:5}

Numerical calculation suggests that we can replace the upper bound of the sum in \eqref{eq:dp+p-1} by $p-1$. Namely, the following variation of \eqref{eq:dp+p-1} should be true.
\begin{conj}\label{conj:first}
Let $d\geqslant 2$ be an even integer and let $p\equiv d+1\pmod{2d}$ be a prime. Then
\begin{align*}
\sum_{k=0}^{p-1} (3dk+1)\frac{(\frac{1}{d})_{k}^{2}(\frac{d-1}{d})_{k}(\frac{1}{2d})_{k}}
{k!^{3}4^{k} (\frac{d+2}{2d})_{k}} \equiv 0 \pmod {p^3}.
\end{align*}
\end{conj}

Furthermore, it seems that \eqref{eq:main} is also true modulo $\Phi_n(q)^3$ for $N=n-1$, which we state as the following conjecture
(which is also a generalization of Conjecture \ref{conj:first}).
\begin{conj}\label{conj:second}
Let $d\geqslant 2$ be an even integer and $e$ an indeterminate. Let $n\equiv d+1\pmod {2d}$ be a positive integer. Then
\begin{align}
\sum_{k=0}^{n-1}
[3dk+1] \frac{(q,e,q^{1+d}/e; q^{2d})_{k}(q,q,q^{d-1};q^{d})_{k}}
{(q^{d},e,q^{1+d}/e;q^{d})_{k}(q^{2d},q^{2d},q^{d+2};q^{2d})_{k}}q^{dk}\equiv 0 \pmod{\Phi_{n}(q)^3}.  \label{eq:conj2}
\end{align}
\end{conj}

It should be pointed out that the $d=2$ case of \eqref{eq:conj2} was already proved by Liu and Wang themselves \cite[Theorem 1]{LiuW3}.
For $d=4$, since $(q^{d-1};q^d)_k/(q^{d+2};q^{2d})_k=1/(-q^3;q^4)_k$, one can easily see that each $k$-th summand on the left-hand side of \eqref{eq:main}
is congruent to $0$ modulo $\Phi_n(q)^3$ for $(3n-1)/4<k\leqslant n-1$. Therefore, By Theorem \ref{thm:main}, the $q$-supercongruence \eqref{eq:conj2}
is also true for $d=4$. However, the same arguments do not work for $d\geqslant 6$.

We find that Theorem \ref{thm:main2} for $d=4$ can be further strengthened as follows.
\begin{conj}
Let $n\equiv 5\pmod{8}$ be a positive integer and $e$ an indeterminate. Then
\begin{align*}
\sum_{k=0}^{(3n+1)/8}[12k-1]\frac{(q^{-1},e,q^{3}/e; q^{8})_{k}(q;q^4)_{k}^3 q^{dk} }
{(q^4,e,q^{3}/e;q^4)_{k}(q^6;q^8)_{k}^3 }
\equiv 0 \pmod{\Phi_n(q)^4}.
\end{align*}
In particular, for any prime $p\equiv 5\pmod{8}$,
$$
\sum_{k=0}^{(3p+1)/8}(12k-1)\frac{(-\frac{1}{8})_k (\frac{1}{4})_k^3}{k!4^k(\frac{3}{4})_k^3}\equiv 0\pmod{p^4}.
$$
\end{conj}

Recently, the author and Zudilin \cite{GuoZu2} have extended many classical $q$-supercongruences to the so-called Dwork-type $q$-supercongruences
through a creative $q$-microscope. They also proposed several difficult conjectures on Dwork-type $q$-supercongruences.
Here we would like to propose such extensions of \eqref{eq:6k+1-LW} and \eqref{conj:cases} for $n\equiv 1\pmod{4}$.
We notice that a similar conjecture related to \eqref{eq;qG2} was already made by Liu and Wang \cite{LiuW2}.

\begin{conj}
Let $n\equiv 1\pmod{4}$ be a positive integer and let $r\geqslant 1$. Then, modulo $[n^r]\prod_{j=1}^r\Phi_{n^j}(q)^2$,
\begin{align*}
&\sum_{k=0}^{(n^r-1)/d}[6k+1]\frac{(q;q^{4})_k(q;q^2)_k^3}{(q^2;q^2)_k(q^{4};q^{4})_k^3}q^{k^2+k} \notag \\
&\quad\equiv \frac{(q^2;q^4)_{(n^r-1)/4}(q^{4n};q^{4n})_{(n^{r-1}-1)/4}}
{(q^4;q^4)_{(n^r-1)/4}(q^{2n};q^{4n})_{(n^{r-1}-1)/4}}[n]q^{(1-n)/4}  \\
&\qquad\times\sum_{k=0}^{(n^{r-1}-1)/d}[6k+1]
\frac{(q^{n};q^{4n})_k(q^n;q^{2n})_k^3}{(q^{2n};q^{2n})_k(q^{4n};q^{4n})_k^3}q^{(k^2+k)n},
\end{align*}
where $d=1,2$.
\end{conj}

\begin{conj}
Let $n\equiv 1\pmod{4}$ be a positive integer and let $r\geqslant 1$. Then, modulo $\Phi_{n^r}(q)\prod_{j=1}^r\Phi_{n^j}(q)^2$,
\begin{align*}
&\sum_{k=0}^{(n^r-1)/2}[6k-1]\frac{(q^{-1};q^{4})_k(q^{-1};q^2)_k(q;q^2)_k^2}{(q^2;q^2)_k(q^{4};q^{4})_k(q^{2};q^{4})_k^2}q^{k^2+k+1} \notag \\
&\quad\equiv \frac{(q^2;q^4)_{(n^r-1)/4}(q^{4n};q^{4n})_{(n^{r-1}-1)/4}}
{(q^4;q^4)_{(n^r-1)/4}(q^{2n};q^{4n})_{(n^{r-1}-1)/4}} q^{(n-1)/4}  \\
&\qquad\times\sum_{k=0}^{(n^{r-1}-1)/2}[6k-1]
\frac{(q^{-n};q^{4n})_k(q^{-n};q^{2n})_k(q^n;q^{2n})_k^2}{(q^{2n};q^{2n})_k(q^{4n};q^{4n})_k(q^{2n};q^{4n})_k^2}q^{(k^2+k+1)n}.
\end{align*}
\end{conj}


\end{document}